\documentclass[twoside]{article}

\usepackage[accepted]{aistats2024}

\usepackage[round]{natbib}

\usepackage{hyperref}
\usepackage[hyphenbreaks]{breakurl}
\usepackage{amsmath,amssymb,amsthm,bm}
\usepackage{multirow}
\usepackage{tikz}
\usepackage{graphicx}
\usepackage{algorithm}
\usepackage{algpseudocode}
\usepackage{caption}
\floatname{algorithm}{Algorithm}

\theoremstyle{plain}
\newtheorem{definition}{Definition}
\newtheorem{theorem}{Theorem}
\newtheorem{lemma}{Lemma}
\newtheorem{proposition}{Proposition}
\newtheorem{corollary}{Corollary}
\newtheorem{assumption}{Assumption}
\newtheorem{example}{Example}
\newtheorem{remark}{Remark}

\DeclareMathOperator*{\argmin}{arg\,min}
\def\x{{\mathbf x}}
\def\y{{\mathbf y}}
\def\z{{\mathbf z}}

\def\R{{\mathbb{R}}}

\begin{document}

\twocolumn[

\aistatstitle{Near-Optimal Convex Simple Bilevel Optimization with a Bisection Method}

\aistatsauthor{Jiulin Wang \And Xu Shi \And Rujun Jiang$^{\dag}$}

\aistatsaddress{School of Data Science, Fudan University}]
\begin{abstract}
This paper studies a class of simple bilevel optimization problems where we minimize a composite convex function at the upper-level subject to a composite convex lower-level problem. Existing methods either provide asymptotic guarantees for the upper-level objective or attain slow sublinear convergence rates. We propose a bisection algorithm to find a solution that is $\epsilon_f$-optimal for the upper-level objective and $\epsilon_g$-optimal for the lower-level objective. In each iteration, the binary search narrows the interval by assessing inequality system feasibility. Under mild conditions, the total operation complexity of our method is ${\tilde  {\mathcal{O}}}\left(\max\{\sqrt{L_{f_1}/\epsilon_f},\sqrt{L_{g_1}/\epsilon_g} \} \right)$. Here, a unit operation can be a function evaluation, gradient evaluation, or the invocation of the proximal mapping, $L_{f_1}$ and $L_{g_1}$ are the Lipschitz constants of the upper- and lower-level objectives' smooth components, and  ${\tilde {\mathcal{O}}}$ hides logarithmic terms.  Our approach achieves a near-optimal rate, matching the optimal rate in unconstrained smooth or composite convex optimization when disregarding logarithmic terms. Numerical experiments demonstrate the effectiveness of our method. 
\end{abstract}

\section{INTRODUCTION}
In this paper, we focus on the following convex bilevel optimization problem:
\begin{equation}\label{pb:BP}
\begin{array}{lcl}
({\rm P}) &\min\limits_{\x\in \R^n}& f(\x):=f_1(\x)+f_2(\x)    \\
&{\rm s.t.}&\x\in \argmin\limits_{\z\in \R^n} g (\z):=g_1(\z)+g_2(\z).
\end{array}
\end{equation}

Here, functions $f_1$ and $g_1: X\rightarrow \R$ are convex and continuously differentiable over an open set $X\in\R^n$. Their gradients, $\nabla f_1$ and $\nabla g_1$, are $L_{f_1}$- and $L_{g_1}$-Lipschitz continuous, respectively. $f_2$ and $g_2: \R^n\rightarrow \R\cup {\infty}$ are proper lower semicontinuous (l.s.c.) convex functions. We assume that $g$ is not strongly convex, and the lower-level problem has multiple optimal solutions; in other words, the optimal solution set of the lower-level problem, denoted as $X_{g}^*$, is not a singleton. Otherwise, the optimal minimum is determined by the lower-level problem. 

This specific class of problems, often known as ``simple bilevel optimization’’ in the existing literature \citep{dempe2010optimality, dutta2020algorithms,shehu2021inertial,jiang2022conditional}, is a subclass of the general bilevel optimization problems. In a general bilevel optimization, the lower-level problem is parametrized by some upper-level variables. Bilevel optimization has garnered significant interest owing to its versatile applications across domains such as reinforcement learning \citep{hong2020two}, meta-learning \citep{bertinetto2018meta, rajeswaran2019meta}, hyper-parameter optimization \citep{franceschi2018bilevel, shaban2019truncated}, and adversarial learning \citep{bishop2020optimal, wang2021fast, wang2022solving}.

Let $p^*$ be the optimal value of problem (\ref{pb:BP}) and $g^*$ be the optimal value of the unconstrained lower-level problem
\begin{equation}\label{pb:g}
\min_{\x\in \R^n}~ g(\x):=g_1(\x)+g_2(\x).
\end{equation}
The goal of this paper is to find an $(\epsilon_f, \epsilon_g)$-optimal solution $\hat{\x}$ satisfying
\begin{equation*}
f(\hat{\x}) - p^* \leq \epsilon_f \quad \text{and} \quad g(\hat{\x}) - g^* \leq \epsilon_g.
\end{equation*}
A possible approach for solving problem (\ref{pb:BP}) is to reformulate it to a constrained optimization problem with functional constraints and apply primal-dual methods.
Specifically, problem (\ref{pb:BP}) can be reformulated as a constrained convex optimization problem as follows:
\begin{equation}
\label{pb:BP1}
\min_{\x\in \R^n} f(\x)~~{\rm s.t.}~~g (\x)\le g^*.
\end{equation}
A critical issue of applying primal-dual-type methods is that problem (\ref{pb:BP1}) does not satisfy the regularity condition required for their convergence (the strict feasibility does not hold and hence Slater's condition fails). Furthermore, classical first-order algorithms, such as projected gradient descent, may also be ineffective due to the difficulty of computing the orthogonal projection onto the level-set of the lower-level objective. If we relax the constraint and solve the following problem to ensure strict feasibility
\begin{equation}
\label{pb:BP2}
\min_{\x\in \R^n} f(\x)~~{\rm s.t.}~~g (\x)\le g^*+\epsilon,
\end{equation}
these challenges remain. Indeed, as $\epsilon$ approaches zero, causing the problem to become nearly degenerate, the dual optimal variable may tend towards infinity. This phenomenon hinders convergence and results in numerical instability \citep{bonnans2013perturbation}. Consequently, problem (\ref{pb:BP}) cannot be straightforwardly addressed as a conventional constrained optimization problem; instead, it necessitates novel theories and algorithms customized for its hierarchical structure.

\subsection{Our Approach}\label{ourapproach}
Our main technique is a bisection method that iteratively narrows an interval $[l,u]$ that includes $p^*$.
The binary search is based on the feasibility of the following system:
\begin{equation}\label{system1}
f(\x)\le c,~ g(\x)\le g^*,
\end{equation}
where $c=\frac{l+u}{2}$.
The key observation is that if System (\ref{system1}) is feasible, then $c$ is an upper bound of $p^*$;
otherwise, System (\ref{system1}) is infeasible and $c$ is a lower bound of $p^*$. This process divides the interval in half.
The feasibility of the system can be checked by solving the following problem
\begin{equation}
\label{pb:system}
\min_{\x\in \R^n}~g(\x),~~{\rm s.t.}~~ f(\x)\leq c.
\end{equation}
The above is only a basic idea, the detail of our algorithm that considers the inexactness of solving (\ref{pb:system}) is detailed in Section \ref{algorithm}. Moreover, by showing that each iteration and the initial lower and upper bounds can be solved by Accelerated Proximal Gradient (APG) methods \citep{nesterov1983method,Beck2009}, we derive a comprehensive complexity analysis for our algorithm.

We state our contributions in the following:
\begin{itemize}
\item Under mild conditions, we propose a novel bisection method that finds an $(\epsilon_f,\epsilon_g)$-optimal solution of problem (\ref{pb:BP}) with an operation complexity $\tilde{\mathcal{O}}\left(\max\{\sqrt{L_{f_1}/\epsilon_f},\sqrt{L_{g_1}/\epsilon_g}\} \right)$, where the notation $\tilde {\mathcal{O}}$ suppresses a logarithmic term. 
Our method achieves near-optimal non-asymptotic guarantees on both upper- and lower-level problems, i.e., our rate aligns with the optimal rate observed in unconstrained smooth or composite convex optimization, with the exception of omitting the logarithmic term \citep{nemirovskij1983problem,woodworth2016tight}.
\item With an additional $r$-th-order ($r\ge 1$) H{\"{o}}lderian error bound assumption on the lower-level problem and incorporating other smoothness assumptions, our method can find a solution $\hat{\x}$ that satisfies $|f(\hat{\x})-p^*|\le \epsilon_f$ and $g(\hat{\x})-g^*\le \epsilon_g$ with an $\tilde{\mathcal{O}}\left(1/\sqrt{\epsilon_f^r} \right)$ operation complexity. This complexity arises under the setting $\epsilon_g=\frac{\alpha}{\gamma} \left(\frac{\epsilon_f}{B_f}\right)^r
$, where $\alpha, \gamma, B_f$ are defined in Section \ref{sec3.4}.
\item  With an additional assumption on the optimal value of (\ref{pb:BP2}), our method can find a solution $\hat{\x}$ that satisfies $|f(\hat{\x})-p^*|\le \epsilon_f$ and $g(\hat{\x})-g^*\le \epsilon_g$ with an $\tilde {\mathcal{O}}\left( \max \{1/\sqrt{L_{\epsilon_g}\epsilon_g}, 1/\sqrt{\epsilon_g}\}\right)$ operation complexity. This operational complexity is observed when $\epsilon_f= L_{\epsilon_g}\epsilon_g$, where $L_{\epsilon_g}$ is defined in Section \ref{sec3.5}.
\item Numerical experiment results on different problems demonstrate the superior performance of our method compared to the state-of-the-art.
\end{itemize}

\begin{table*}[ht]
\centering
\caption{Summary of simple bilevel optimization algorithms. The abbreviations ``SC'', ``C'', ``C3''
stand for ``strongly convex'', ``convex'', ``Convex objective with Convex Compact constraints''  respectively. When the connection between complexity and the gradient's Lipschitz constant is clear, we include it in the complexity result; otherwise, we omit it.}
\label{table1}
\resizebox{1.0\textwidth}{!}{
\begin{tabular}{cccccc}
\hline
\multirow{2}{*}{References} & Upper-level &{Lower-level} & \multicolumn{2}{c}{Convergence} \\
\cline{2-5}
& Objective $f$  &Objective $g$  & Upper-level & Lower-level\\
 \hline
MNG \citep{beck2014first}.
& SC, differentiable  & C, smooth   & Asymptotic  & ${\mathcal{O}}\left( L_{g_1}^2/\epsilon^2\right)$ \\
 \hline
BiG-SAM \citep{sabach2017first}&SC, smooth  &C, composite    &Asymptotic  &${\mathcal{O}}\left( L_{g_1}/\epsilon\right)$ \\
 \hline
IR-IG \citep{amini2019iterative}&SC  &C3, Finite
sum   &Asymptotic  &${\mathcal{O}}\left( 1/\epsilon^4\right)$ \\
 \hline
Tseng's method \citep{malitsky2017chambolle}&C, composite  &C, composite    &Asymptotic  &${\mathcal{O}}\left( 1/\epsilon\right)$ \\
\hline
ITALEX \citep{doron2022methodology}&C, composite  &C, composite     &${\mathcal{O}}\left(1/\epsilon^2\right)$  &${\mathcal{O}}\left(1/\epsilon\right)$  \\
 \hline
a-IRG \citep{kaushik2021method}&C, Lipschitz  &C, Lipschitz    &\multicolumn{2}{c}{${\mathcal{O}}\left(\max\{1/\epsilon_f^4,1/\epsilon_g^4\} \right)$ } \\
\hline
CG-BiO \citep{jiang2022conditional}&C, smooth  &C3, smooth & \multicolumn{2}{c}{${\mathcal{O}}\left(\max\{L_{f_1}/\epsilon_f,L_{g_1}/\epsilon_g\} \right)$ }  \\
 \hline
{\bf Our method} &C, composite  &C, composite   &\multicolumn{2}{c}{ $\tilde{\mathcal{O}}\left(\max\{\sqrt{L_{f_1}/\epsilon_f},\sqrt{L_{g_1}/\epsilon_g}\} \right)$ }  \\
 \hline
 \end{tabular}
}
\end{table*}

\subsection{Related Work}
One class of algorithms to solve problem (\ref{pb:BP}) is based on solving the Tikhonov-type regularization \citep{Tikhonov1977}:
\begin{equation*}
\min_{\x\in \R^n} \phi(\x):=g (\x)+\lambda f(\x),
\end{equation*}
where $\lambda>0$ is a regularization parameter. However, these kinds of algorithms fail to provide any non-asymptotic guarantee for either the upper- or lower-level objective. For a review of these algorithms, see \cite{doron2022methodology} and \cite{jiang2022conditional}.

Another class of algorithms aims to establish non-asymptotic convergence rates for problem (\ref{pb:BP}). \cite{beck2014first} presented the Minimal Norm Gradient (MNG) method for the case where $f$ is strongly convex. They demonstrated that MNG converges asymptotically to the optimal solution and possesses an ${\mathcal{O}}\left({L_{g_1}^2}/{\epsilon^2} \right)$ complexity bound for the lower-level problem.  In their setting, $g\equiv g_1$ and $g_2\equiv 0$. Developed from the sequential averaging method (SAM) framework, the Bilevel Gradient Sequential Averaging Method (BiG-SAM) is proposed by \citet{sabach2017first}. This algorithm can achieve an ${\mathcal{O}}\left( L_{g_1}/\epsilon\right)$ complexity bound for the lower-level problem. \citet{solodov2007explicit} introduced the Iterative Regularized Projected Gradient (IR-PG) method, which involves applying a projected gradient step to the Tikhonov-type regularization function $\phi(\x)$
at each iteration. \citet{amini2019iterative} extended the IR-PG method \citep{solodov2007explicit} for the case where $f$ is strongly convex but not necessarily differentiable. Their method achieves a convergence rate of ${\mathcal{O}}\left( 1/k^{0.5-b}\right)$ for the lower-level problem, where $b\in (0,0.5)$. \citet{malitsky2017chambolle} studied a version of Tseng's accelerated gradient method that obtains a convergence rate of ${\mathcal{O}}\left( 1/k\right)$ for the lower-level problem. These prior works only establish the convergence rate for the lower-level problem, while the rate for the upper-level objective is missing.

Several algorithms have recently provided convergence rates for both upper- and lower-level objectives. \citet{doron2022methodology} presented a scheme called Iterative Approximation and Level-set Expansion (ITALEX) to solve problem (\ref{pb:BP}). Their algorithm achieves convergence rates of ${\mathcal{O}}\left( 1/k\right)$ and ${\mathcal{O}}\left( 1/\sqrt{k}\right)$ for the lower- and upper-level problems, respectively.  \citet{kaushik2021method} showed that an iteratively regularized gradient (a-IRG) method can obtain complexity  ${\mathcal{O}}\left( 1/k^{0.5-b}\right)$ for the upper-level problem and ${\mathcal{O}}\left( 1/k^{b}\right)$ for the lower-level, where $b\in (0,0.5)$. To balance the two rates, one can set $b=0.25$, and the complexity bound is ${\mathcal{O}}\left(\max\{1/\epsilon_f^4,1/\epsilon_g^4\} \right)$ as stated in Table \ref{table1}. \citet{jiang2022conditional} presented a conditional gradient-based bilevel optimization (CG-BiO) method, which requires ${\mathcal{O}}\left(\max\{L_{f_1}/\epsilon_f, L_{g_1}/\epsilon_g\} \right)$ operation complexity to find an $(\epsilon_f,\epsilon_g)$-optimal solution.  In their setting, $f\equiv f_1$ and $f_2\equiv 0$.

\section{PRELIMINARIES}\label{preliminaries}
In this paper, we use an Accelerated Proximal Gradient (APG) method to approximately solve subproblems, which have the following form:
\begin{equation}\label{pb:varphi}
\min_{\x\in \R^n}~ \varphi(\x):=\varphi_1(\x)+\varphi_2(\x),
\end{equation}
where the function $\varphi_1: X\rightarrow \R$ is convex and continuously differentiable on an open set $X\in\R^n$. The  gradient $\nabla \varphi_1$ is $L_{\varphi_1}$-Lipschitz continuous. The function $\varphi_2: \R^n\rightarrow \R\cup \{\infty\}$ is proper, lower semicontinuous, convex, possibly non-smooth, and proximal-friendly. A function $h$ is proximal-friendly means that the proximal mapping of $h$, defined as
\[
{\rm prox}_h(\y)=\arg \min_{\x\in \R^n} h(\x)+\frac{1}{2}\|\x-\y\|^2,
\]
is easy to compute. In this paper, we take the classical Fast Iterative Shrinkage
Thresholding Algorithm (FISTA) proposed in \citet{Beck2009} as an APG algorithm (see more details in the appendix). Next, we give a definition of an APG oracle.
\begin{definition}
\label{definition-APG}
Given $\varphi_1$, $\varphi_2$, $L_{\varphi_1}$, and $\x_0^{\varphi} \in \R^n$ as defined above, an APG oracle, denoted by $\tilde{\x}^{\varphi} = \text{APG}(\varphi_1, \varphi_2, L_{\varphi_1}, \x_0^{\varphi}, \epsilon)$, is a procedure that implements the classical FISTA scheme within $\mathcal{O}(\sqrt{L_{\varphi_1}/\epsilon})$ iterations to obtain an $\epsilon$-optimal solution to problem (\ref{pb:varphi}), denoted as $\tilde{\x}_{\varphi}$.
\end{definition}
If the Lipschitz constant $L_{\varphi_1}$ is unknown or computationally infeasible, we can apply the FISTA scheme with line search as an alternative to the APG algorithm  (see the appendix).

For any fixed $c$, problem (\ref{pb:system}) can be rewritten in the following form:
\begin{equation}
\label{pb:c}
\min_{\x\in \R^n}~ g_c(\x):=g_1(\x)+h_c(\x),
\end{equation}
where $h_c(\x):=g_2(\x)+\delta_c(\x)$ and $\delta_c(\x)$ is an indicator function with the definition that $\delta_c(\x)=0$ if $f(\x)\le c$; $\delta_c(\x)=+\infty$ if $f(\x)> c$. We provide some motivating examples in which the proximal mapping of $h_c$ is easy to compute. These examples can be found in the appendix.

\subsection{Assumptions}
\begin{assumption}\label{as1}
We adopt the following basic assumptions.
\begin{itemize}
\item [(i)] Functions $f_1$ and $g_1$ are convex and continuously differentiable. The gradients of the functions $f_1$, $g_1$, denoted by $\nabla {f_1}$ and $\nabla {g_1}$ are $L_{f_1}$- and $L_{g_1}$-Lipschitz continuous, respectively.
\item[(ii)]  Functions $f_2$ and $g_2$ are proper, lower semicontinuous, convex, possibly non-smooth, and proximal-friendly.
\item[(iii)] The upper- and lower-level functions are lower bounded:
\[
f^*:=\inf_{\x\in \R^n} f(\x)>-\infty,~g^*:=\inf_{\x\in \R^n} g(\x)>-\infty.
\]
In addition, we assume that $g$ is not strongly convex, and the lower-level problem has multiple optimal solutions.
\item[(iv)] For any fixed $c$, the function $h_c:=g_2+\delta_c$ in problem (\ref{pb:c}) is proximal-friendly.
\end{itemize}
\end{assumption}
\begin{remark}
We assume that the upper-level problem involves the minimization of a composite convex function, which comprises the sum of a smooth convex function and a potentially non-smooth convex function. This assumption is significantly weaker than the strong-convexity assumption made in certain previous studies \citep{beck2014first,sabach2017first,amini2019iterative}. This assumption is also less restrictive than the requirement for the upper-level objective function to be smooth \citep{jiang2022conditional}. We also assume that the lower-level problem is a composite convex minimization. This assumption is less restrictive than the smoothness assumption made in \citet{beck2014first}. Additionally, this assumption is weaker than the requirement that the lower-level objective function is convex with convex compact constraints, as described in \citet{amini2019iterative} and \citet{jiang2022conditional}.
\end{remark}
\begin{remark}
We make Assumption \ref{as1}(iv) to enable the efficient application of the APG oracle for solving problem (\ref{pb:system}). The function $h_c$ is the sum of two convex functions, and the study of proximal mapping for such sums can be found in the literature \citep{yu2013decomposing,pustelnik2017proximity,bauschke2018projecting,adly2019decomposition}.
\end{remark}

\section{MAIN ALGORITHM AND CONVERGENCE ANALYSIS}\label{algorithm}
Before the presentation of the main algorithm, we describe our idea in the next subsection.
\subsection{Bisection Method}\label{sec3.0}
\begin{figure}
\tikzset{every picture/.style={line width=0.75pt}}
\begin{tikzpicture}[x=0.75pt,y=0.75pt,yscale=-1,xscale=1]
\draw [color={rgb, 255:red, 0; green, 0; blue, 0 }  ,draw opacity=1 ][line width=0.75]  (41,205.15) -- (280.67,205.15)(61.06,16) -- (61.06,231) (273.67,200.15) -- (280.67,205.15) -- (273.67,210.15) (56.06,23) -- (61.06,16) -- (66.06,23)  ;
\draw    (81.67,35) .. controls (83.67,171) and (128.67,171) .. (171,171) ;
\draw    (171,171) -- (240.67,171) ;
\draw  [dash pattern={on 0.84pt off 2.51pt}]  (62,151) -- (108,151) ;
\draw  [dash pattern={on 0.84pt off 2.51pt}]  (62,171) -- (171,171) ;
\draw  [dash pattern={on 0.84pt off 2.51pt}]  (108,151) -- (108,206) ;
\draw  [dash pattern={on 0.84pt off 2.51pt}]  (78,25.5) -- (78,206) ;
\draw  [dash pattern={on 0.84pt off 2.51pt}]  (171,171) -- (171,206) ;
\draw  [dash pattern={on 0.84pt off 2.51pt}]  (206,171) -- (206,206) ;
\draw  [dash pattern={on 0.84pt off 2.51pt}]  (108,151) -- (171,171) ;
\draw (108,151) circle (3);
\draw (171,171) circle (3);
\draw (206,171) circle (3);
\draw (92,49) node [anchor=north west][inner sep=0.75pt]   [align=left] {$\displaystyle \overline{g}( c)$};
\draw (68,207) node [anchor=north west][inner sep=0.75pt]   [align=left] {$\displaystyle f^{*}$};
\draw (99,207) node [anchor=north west][inner sep=0.75pt]   [align=left] {$\displaystyle p_{\epsilon_g}^{*}$};
\draw (163,207) node [anchor=north west][inner sep=0.75pt]   [align=left] {$\displaystyle p^{*}$};
\draw (184,207) node [anchor=north west][inner sep=0.75pt]   [align=left] {$\displaystyle p^{*} +\epsilon _{f}$};
\draw (14,143) node [anchor=north west][inner sep=0.75pt]   [align=left] {$\displaystyle g^{*}+\epsilon _{g}$};
\draw (41,161) node [anchor=north west][inner sep=0.75pt]   [align=left] {$\displaystyle g^{*}$};
\draw (269,210) node [anchor=north west][inner sep=0.75pt]   [align=left] {$\displaystyle c$};
\end{tikzpicture}
\caption{Variation of $\bar{g}(c)$ over $(f^{*},+\infty)$}
\label{figure1}
\end{figure}
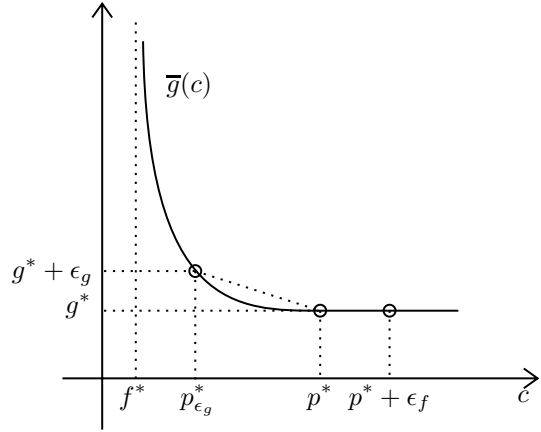

Our method is a bisection method whose heart is the feasibility of System (\ref{system1}). Let $f^*$ be the optimal value of the unconstrained upper-level problem
\begin{equation}\label{pb:fvalue}
\min_{\x\in \R^n}~ f(\x):=f_1(\x)+f_2(\x).
\end{equation}
For a given $c>f^*$, we let $\bar{g} (c)$ be the optimal value of problem (\ref{pb:system}). Then $\bar{g} (c)$ is a univariate function of $c$ on $(f^*,+\infty)$. According to Theorem 5.3 in \citet{Rockafellar1970}, the function $\bar{g} (c)$ is convex. The function $\bar{g} (c)$ is also non-increasing as the feasible set of problem (\ref{pb:system}) becomes larger when $c$ increases. Moreover, if $f^*<c<p^*$, then the inequality $\bar{g} (c)>g^*$ holds; otherwise $c\ge p^*$ and we have $\bar{g} (c)=g^*$. Therefore, $p^*$ is the left-most root of the equation $\bar{g}(c)=g^*$. Let $p_{\epsilon_g}^*$ be the optimal value of (\ref{pb:BP2}) with $\epsilon=\epsilon_g$, then it is a root of the equation $\bar{g}(c)=g^*+\epsilon_g$. We illustrate the graph of $\bar{g} (c)$ in Figure \ref{figure1}.

To illustrate the basic idea of our method, we make an ideal assumption that the exact values of $g^*$ and $\bar{g}(c)$ can be obtained. We can observe that if $\bar{g} (c)>g^*$ holds, then System (\ref{system1}) is infeasible; otherwise, $\bar{g} (c)=g^*$ and System (\ref{system1}) is feasible. For a guess point $c$, if the condition $\bar{g} (c)>g^*$ holds, then $c$ is a lower bound of $p^*$; otherwise, $c$ is an upper bound of $p^*$.

However, the ideal assumption that the exact values of $g^*$ and $\bar{g} (c)$ can be obtained does not hold. Instead, we solve problem (\ref{pb:g}) and problem (\ref{pb:system}) to approximate them, respectively. For problem (\ref{pb:g}), we invoke the APG oracle ${\tilde{\x}_g}$=APG$(g_1, g_2, L_{g_1}, \x_0^g, \epsilon_g/2)$ to solve it. Then
we can find an approximate solution $\tilde{\x}_{g}$ that satisfies
\begin{equation}\label{approxg}
0\le \tilde{g}-g^*\leq \epsilon_g/2,
\end{equation}
where $\tilde{g}=g(\tilde{\x}_{g})$.

Under Assumption \ref{as1}(iv), the proximal mapping of $h_c$ is easy to compute. Then we can apply an APG oracle to problem (\ref{pb:system}). For a given $c$, we invoke the APG oracle ${\tilde{\x}_c}$=APG$(g_1, h_c, L_{g_1}, \x_0^c, \epsilon_g/2)$ to solve problem (\ref{pb:system}). Then we can obtain an approximate solution $\tilde{\x}_{c}$ that satisfies
\begin{equation}\label{approxpsi}
0\le g(\tilde{\x}_c)-\bar{g} (c)\leq \epsilon_g/2.
\end{equation}
Since the condition $\bar{g} (c)>g^*$ cannot be verified directly, we replace it with the following verifiable condition
\begin{equation}\label{cond1}
g(\tilde{\x}_c)> \tilde{g}+\epsilon_g/2.
\end{equation}
If Condition (\ref{cond1}) holds, then we have
\[
\bar{g}(c)\overset{(\ref{approxpsi})}{\ge } g(\tilde{\x}_c)-\epsilon_g/2\overset{(\ref{cond1})}{>}\tilde{g}\overset{(\ref{approxg})}{\ge} g^*,
\]
i.e., the inequality $\bar{g}(c)>g^*$ holds. Thus, System (\ref{system1}) is infeasible, and $c$ is a lower bound of $p^*$.

If Condition (\ref{cond1}) does not hold, we have $g(\tilde{\x}_c)\le \tilde{g}+\epsilon_g/2$ and thus
\begin{equation}\label{approxgbar}
\bar{g}(c)\overset{(\ref{approxpsi})}{\le} g(\tilde{\x}_c)\le \tilde{g}+\epsilon_g/2 \overset{(\ref{approxg})}{\le } g^*+\epsilon_g.
\end{equation}
We cannot infer that $\bar{g}(c)=g^*$ holds in this case. Also, we cannot claim that System (\ref{system1}) is feasible. Thus $c$ might not be an upper bound of $p^*$.  By (\ref{approxgbar}), $\tilde{\x}_c$ is a feasible solution of (\ref{pb:BP2}) with $\epsilon=\epsilon_g$. Therefore, $f(\tilde{\x}_c)$ is an upper bound on $p_{\epsilon_g}^*$, where $p_{\epsilon_g}^*$ represents the optimal value of (\ref{pb:BP2}) with $\epsilon=\epsilon_g$ as previously defined. In addition, as the inequality $ g(\tilde{\x}_c)\le g^*+\epsilon_g$ holds, $\tilde{\x}_c$ is an $\epsilon_g$-optimal solution of the lower-level problem (\ref{pb:g}).

Summarizing the above analysis, we present the following lemma.

\begin{lemma}\label{lem:system}
For any fixed $c$, if Condition (\ref{cond1}) is satisfied, then System (\ref{system1}) is infeasible, and $c$ is a lower bound of $p^*$. If Condition (\ref{cond1}) is not satisfied, then we can obtain $\tilde{\x}_c$ as an $\epsilon_g$-optimal solution of the lower-level problem and $f(\tilde{\x}_c)$ is an upper bound of $p_{\epsilon_g}^*$.
\end{lemma}

Next, we show how to obtain the initial interval $[l,u]$ for the bisection procedure. Here $l$ is a lower bound of $p^*$ and $u$ is an upper bound of $p_{\epsilon_g}^*$, but possibly not an upper bound of $p^*$. We invoke the APG oracle ${\tilde{\x}_f}$=APG$(f_1, f_2, L_{f_1}, \x_0^f, \epsilon_f/2)$ to solve problem (\ref{pb:fvalue}). Then we can obtain an approximate solution $\tilde{\x}_{f}$ that satisfies
\begin{equation}\label{approxfvalue}
0\le f(\tilde{\x}_{f})-f^*\leq \epsilon_f/2.
\end{equation}
We have $f(\tilde{\x}_{f})-\epsilon_f/2\le f^*\le p^*$. We use $l=f(\tilde{\x}_{f})-\epsilon_f/2$ as an initial lower bound of $p^*$. Moreover, since $\tilde{\x}_{g}$ satisfies (\ref{approxg}), we have $\tilde{\x}_{g}$ is a feasible solution of (\ref{pb:BP2}) with $\epsilon=\epsilon_g$. Then we use $u=f(\tilde{\x}_{g})$ as an initial upper bound of $p_{\epsilon_g}^*$.

Now we can do a binary search over $[l,u]$. For a given $c=\frac{l+u}{2}$, we check whether Condition (\ref{cond1}) is satisfied.
If Condition (\ref{cond1}) is satisfied, we let $l=c$ be a new lower bound of $p^*$. If Condition (\ref{cond1}) is not satisfied, we let $u=f(\tilde{\x}_c)$, which is less than or equal to $c$, be a new upper bound of $p_{\epsilon_g}^*$.
We summarise our method in Algorithm \ref{main algorithm}.

\begin{algorithm}
\caption{Bisection-based method for simple Bilevel Optimization (Bisec-BiO)}
\label{main algorithm}
\begin{algorithmic}[1]
\Require $f_1,f_2,g_1,g_2,L_{f_1},L_{g_1},\epsilon_f,\epsilon_g$
\Ensure An ($\epsilon_f,\epsilon_g$)-optimal solution $\hat{\x}$

\State Invoke APG oracles to obtain initial bounds $l$ and $u$, and the approximate solutions $\tilde{\x}_{f}$ and $\tilde{\x}_{g}$.
\While{\parbox[t]{0.20\linewidth}{$u-l>\epsilon_f$}}
        \State  let $c=\frac{l+u}{2}$ and invoke an APG oracle to obtain an approximate solution $\tilde{\x}_c$.
        \If{\parbox[t]{0.52\linewidth}{Condition (\ref{cond1}) is satisfied}}
            \State let $l=c$,
        \Else
            \State let $u=f(\tilde{\x}_c)$.\qquad\qquad\qquad$\triangleright~ f(\tilde{\x}_c)\le c$
        \EndIf
\EndWhile

\State Let $c=u$ and return
the corresponding $\tilde{\x}_c$ as $\hat{\x}$.
\end{algorithmic}
\end{algorithm}

\subsection{Convergence Analysis under Assumption \ref{as1}}\label{sec3.3}
In this subsection, we give the complexity result of our method.
\begin{theorem}\label{thm1}
Suppose Assumption \ref{as1} holds. Algorithm \ref{main algorithm} produces an
($\epsilon_f,\epsilon_g$)-optimal solution for problem (\ref{pb:BP}) after at most $T$ evaluations of the function values $f_1$, $f_2$, $g_1$ and $g_2$, the gradients $\nabla f_1$ and $\nabla g_1$, and the calls of proximal mapping with respect to function $h_c$, where
\[
T=\tilde {\mathcal{O}}\left(\max\left\{\sqrt{\frac{L_{f_1}}{ \epsilon_f}},\sqrt{\frac{L_{g_1}}{ \epsilon_g}} \right\}\right),
\]
and $\tilde {\mathcal{O}}$ suppresses a logarithmic term.
\end{theorem}
Theorem \ref{thm1} demonstrates that our complexity achieves the near-optimal rate for both upper- and lower-level objectives and matches the optimal rate of first-order methods for unconstrained smooth or composite convex optimization when disregarding the logarithmic term \citep{nemirovskij1983problem, woodworth2016tight}. Comparing with existing works \citep{beck2014first, sabach2017first, amini2019iterative, malitsky2017chambolle, doron2022methodology, kaushik2021method, jiang2022conditional}, our result provides the best-known non-asymptotic bounds for both upper- and lower-level objectives. Specifically, our complexity bound improves upon the result by \cite{jiang2022conditional} by orders of magnitude. They considered a different setup where the upper-level function is smooth, and the lower-level objective is a smooth convex function with convex compact constraints.

\begin{remark}  
Assumption 1(iv) essentially implies that projecting onto the sublevel set of $f$ is straightforward, as observed in the motivating examples. Consequently, we can employ a bisection method to locate $f^*$ by verifying the solvability of the projection onto the sublevel set of $f$. This approach eliminates the dependency on $L_{f_1}$ from the overall complexity.
\end{remark}
\section{CONVERGENCE ANALYSIS UNDER OTHER ASSUMPTIONS}
In this section, we provide an analysis of the metric $f(\hat{\x}) - p^*$. Our method guarantees that $f(\hat{\x})$ serves as an upper bound for $p_{\epsilon_g}^*$; however, it may not necessarily be an upper bound for $p^*$, which could result in a negative value for $f(\hat{\x}) - p^*$. In fact, as $\hat{\x}$ may not be an exact optimal solution to the lower-level problem, it may not be a feasible point for problem (\ref{pb:BP}). Hence, the value of $f(\hat{\x}) - p^*$ may be negative. In essence, we currently offer an upper bound metric for $f(\hat{\x})-p^*$ while lacking a corresponding lower bound metric. Although this may appear a bit unconventional, it is noteworthy that \citet{kaushik2021method} and \citet{jiang2022conditional} similarly employed $f(\hat{\x})-p^*$ as their performance metric.

In this section, we introduce additional assumptions to establish a lower bound for $f(\hat{\x})-p^*$, allowing us to provide a metric expressed as $|f(\hat{\x})-p^*|$.
\subsection{Convergence Analysis under H{\"{o}}lderian Error Bound Assumption}\label{sec3.4}

In this subsection, we make some additional assumptions to provide a lower bound for $f(\hat{\x})-p^*$.
\begin{assumption}\label{as2}
\begin{itemize}
    \item[(i)] The domain of $g_2$ is bounded.
    \item[(ii)] The function $f_2$ is $l_{f_2}$-Lipschitz continuous on ${\rm dom}(g_2)$, i.e. $|f_2(\x)-f_2(\y)|\le l_{f_2}\|\x-\y\|$.
\end{itemize}
\end{assumption}
In the following, we give some remarks on Assumption \ref{as2}.  Assumption \ref{as2}(i) is fulfilled by many examples, see Section 5 in \citet{amini2019iterative} and Section 2 in \citet{jiang2022conditional}. In particular, in Section \ref{logistic regression problem}, we have $C=\{\x\in \R^n: \|\x\|_1\le \lambda\}$, then ${\rm dom}(g_2)$ is bounded. As $f_1$ is continuous, $\nabla f_1$ is bounded on ${\rm dom}(g_2)$. Hence Assumption \ref{as2}(i) implies
\[
B_{f_1}=\max_{\x\in {\rm dom}(g_2)}  \|\nabla f_1(\x)\|.
\]
Then by mean-value theorem, the function $f_1$ is $B_{f_1}$-Lipschitz continuous on ${\rm dom}(g_2)$.

Assumption \ref{as2}(ii) is mild. For example, if $f_2(\x)=\|\x\|_1$, then $l_{f_2}=\sqrt{n}$, where $n$ is the dimension of $\x$. Let $B_f=B_{f_1}+l_{f_2}$. Under Assumption \ref{as2}(i) and (ii), it follows that $f$ is $B_f$-Lipschitz continuous on ${\rm dom}(g_2)$, namely,
\begin{equation}\label{Lip}
|f(\x_1)-f(\x_2)|\le B_f \|\x_1-\x_2\|~\forall \x_1,\x_2\in {\rm dom}(g_2)
\end{equation}

\begin{assumption}\label{as3}
The function $g$ satisfies the H{\"{o}}lderian error bound for some $\alpha> 0$ and $r\ge1$ on the lower-level optimal solution set $X_g^*$, i.e,
\begin{equation}\label{holderian}
\frac{\alpha}{r}~{\rm dist}(\x,X_g^*)^r\le g(\x)-g^*,~ \forall \x\in {\rm dom}(g_2).  \end{equation}
\end{assumption}
It is important to highlight that the error bound condition described in (\ref{holderian}) has received considerable attention in the literature, as evidenced by studies such as \citet{1997Error, J2017From,zhou2017}, and the associated references therein. \citet{J2017From} demonstrated that this error bound condition typically holds when the function $g$ exhibits properties of being semi-algebraic and continuous, while also ensuring that ${\rm dom}(g_2)$ remains bounded. They also showed that there is an equivalence between H{\"{o}}lderian error bound condition and the Kurdyka-Łojasiewicz inequality. There are two notable special cases: (i) $r = 1$, $X_g^*$ is a set of weak sharp minima of $g$ \citep{2005Weak}; (ii) $r = 2$, Condition (\ref{holderian}) is known as the quadratic growth condition \citep{2018Error}. Based on Corollary 5.1 in \citet{li2018calculus} and Theorem 5 in \citet{J2017From}, it is evident that in the motivating examples provided in the appendix, the lower-level objectives satisfy the H{\"{o}}lderian error bound assumption with $r=2$.

Given Assumptions \ref{as2} and \ref{as3}, we can derive the following lower bound for $f(\hat{\x})-p^*$. Importantly, this result is an inherent characteristic of problem (\ref{pb:BP}) and remains unaffected by the choice of algorithm.
\begin{proposition}\label{prop1}
Under Assumptions \ref{as2} and \ref{as3}, let $\hat{\x}$ be an $\epsilon_g$-optimal solution of the lower-level problem, i.e., $\hat{\x}$ satisfies $g(\hat{\x})-g^*\le \epsilon_g$. Then it holds that:
\[
f(\hat {\x})-p^*\ge -B_f \left(\frac{r\epsilon_g}{\alpha}\right)^{\frac{1}{r}}.
\]
\end{proposition}
This result is similar to Proposition 1 in \cite{jiang2022conditional}, which also gives a lower bound for $f(\hat{\x})-p^*$ with similar assumptions.
Combining Theorem \ref{thm1}  with Proposition \ref{prop1}, we have the following result.
\begin{corollary} \label{cor1}
Under Assumptions \ref{as1}-\ref{as3}, let $\hat {\x}$ be the output of Algorithm \ref{main algorithm} and set $
\epsilon_g=\frac{\alpha}{\gamma} \left(\frac{\epsilon_f}{B_f}\right)^r$. Then with an operation complexity of $\tilde {\mathcal{O}}\left({1/\sqrt{\epsilon_f^r}}\right)$, we can find an $\hat {\x}$ such that
\[
|f(\hat {\x})-p^*|\le \epsilon_f,~ g(\hat {\x})-g^*\le \epsilon_g.
\]
\end{corollary}
Corollary \ref{cor1} illustrates that, under Assumptions \ref{as1}-\ref{as3}, we can find an iteration point to be $\epsilon_f$-close to optimality with an operation complexity $\tilde {\mathcal{O}}\left({1/\sqrt{\epsilon_f^r}}\right)$.
\subsection{Convergence Analysis under an Additional Assumption on the Optimal Value of (\ref{pb:BP2})}\label{sec3.5}
In this subsection, we present an alternative approach to establishing a lower bound for $f(\hat{\x})-p^*$. Denote the optimal value of problem (\ref{pb:BP2}) as $v(\epsilon)$. Here, $v(\epsilon)$ is a function of $\epsilon$ defined on the interval $[0,+\infty)$. According to the definitions of $p^*$ and $p_{\epsilon_g}^*$, we can establish that $v(0)=p^*$ and $v(\epsilon_g)=p_{\epsilon_g}^*$. As illustrated in Figure \ref{figure1}, we can define an angle function $\theta(\epsilon)$ on $[0,\pi/2)$ such that
\[
\epsilon=(v(0)-v(\epsilon))\cdot \tan \theta(\epsilon).
\]
Similar to the monotonicity of $\bar{g}(c)$, functions $v(\epsilon)$ and $\theta(\epsilon)$ are also monotonically non-decreasing. If $\epsilon=0$, then $\theta(\epsilon)=0$; otherwise $\epsilon>0$ and $\theta(\epsilon)>0$.
It is hard to compute the exact value of $\tan \theta(\epsilon)$. Instead, for any fixed $\epsilon>0$, we assume that there exists a parameter $L_{\epsilon}>0$ such that $L_{\epsilon} \tan \theta(\epsilon)\ge 1$ holds. In other words, we make the following assumption.
\begin{assumption}\label{as4}
For any fixed $\epsilon>0$, we assume that there exists a parameter $L_{\epsilon}>0$ such that $v(0)-v(\epsilon)\le L_{\epsilon} \epsilon$ holds.
\end{assumption}
In the following, we consider a simple two-dimensional toy example in which Assumption \ref{as4} holds.
\begin{example}[A toy example]\label{ex3}
\begin{equation}
\label{toy example}
\begin{array}{lcl}
&\min\limits_{\x\in\mathbb{R}^{2}}&f(\x)=|\x_1|+|\x_2|\\
&{\rm{s.t.}}&\x\in\argmin\limits_{\z\in\R^{2}}g(\z) =(\z_1-1)^2.
\end{array}
\end{equation}
\end{example}
In this example, the optimal solution and the optimal value of (\ref{toy example}) are $(\x_1^*,\x_2^*)=(1,0)$ and $v(0)=p^*=1$, respectively. Next, we consider problem (\ref{pb:BP2}) with $0<\epsilon<1$. It can be easily obtained that the optimal solution and the optimal value of (\ref{pb:BP2}) are $(\x_1^*,\x_2^*)=(1-\sqrt{\epsilon},0)$ and $v(\epsilon)=1-\sqrt{\epsilon}$, respectively. Then we have
\[
v(0)-v(\epsilon)=\sqrt{\epsilon}.
\]
By setting $L_\epsilon\ge 1/\sqrt{\epsilon}$, Assumption \ref{as4} holds.

Under Assumption \ref{as4}, we can derive a lower bound for $f(\hat{\x})-p^*$ that is independent of the choice of algorithm.
\begin{proposition}\label{prop2}
Under Assumption \ref{as4}, let $\hat{\x}$ be an $\epsilon_g$-optimal solution of the lower-level problem, i.e., $\hat{\x}$ satisfies $g(\hat{\x})-g^*\le \epsilon_g$. Then it holds that:
\[
f(\hat {\x})-p^*\ge -L_{\epsilon_g} \epsilon_g.
\]
\end{proposition}

By combining Theorem \ref{thm1}  with Proposition \ref{prop2}, we have the following result.
\begin{corollary} \label{cor2}
Under Assumption \ref{as1} and Assumption \ref{as4}, let $\hat {\x}$ be the output of Algorithm \ref{main algorithm} and set $
{\epsilon_f}=L_{\epsilon_g}\epsilon_g$. Then with an operation complexity of $\tilde {\mathcal{O}}\left( \max \{1/\sqrt{L_{\epsilon_g}\epsilon_g}, 1/\sqrt{\epsilon_g}\}\right)$, we can find an $\hat {\x}$ such that
\[
|f(\hat {\x})-p^*|\le \epsilon_f,~ g(\hat {\x})-g^*\le \epsilon_g.
\]
\end{corollary}
Corollary \ref{cor2} demonstrates that with Assumptions \ref{as1} and \ref{as4} in place, we can obtain an iteration point to be $L_{\epsilon_g}\epsilon_g$-close to optimality with an operation complexity $\tilde {\mathcal{O}}\left( \max \{1/\sqrt{L_{\epsilon_g}\epsilon_g}, 1/\sqrt{\epsilon_g}\}\right)$.

\section{NUMERICAL EXPERIMENTS}
\label{experiment}
In this section, we apply our method (Bisec-BiO) to two bilevel optimization problems from the motivating examples in the appendix and compare its performance with other existing methods in the literature \citep{beck2014first,sabach2017first,kaushik2021method,gong2021bi,jiang2022conditional}. For all experiments, we set $\epsilon_f = 10^{-5}$ and $\epsilon_g = 10^{-6}$, and we adopt the Greedy FISTA algorithm proposed in \cite{liang2022improving} as the APG method. The Greedy FISTA algorithm can achieve superior practical performance compared to the classical FISTA.

\subsection{Minimum Norm Solution Problem (MNP)}
\label{minimum norm solution problem}
We first consider the linear regression problem on the YearPredictionMSD dataset\footnote{\url{https://archive.ics.uci.edu/dataset/203/yearpredictionmsd}}, which contains information on $515,345$ songs, with a release year from $1992$ to $2011$. For each song, the dataset contains its release year and an additional $90$ attributes. We use a sample of $1,000$ songs randomly selected from the dataset with uniform i.i.d distribution, and denote the feature matrix and the release years by $A$ and $b$, respectively. Additionally, in line with \citet{merchav2023convex}, we adopt a min-max scaling technique and add an interceptor and $90$ co-linear attributes to $A$.

For the lower-level, we let $g_1(\x) = \frac{1}{2}\|A\x-b\|^2$, which exhibits a $L_{g_1}$-Lipschitz continuous gradient, where $L_{g_1} = \lambda_{\max}(A^TA)$. Simultaneously, we set $g_2(\x) \equiv 0$. For the upper-level, we let $f_1(\x)=\frac{1}{2}\|\x\|^2$ and $f_2(\x)\equiv 0$. This configuration corresponds to finding the minimum norm solution.  Now, our goal is to solve the following bilevel problem:
\begin{equation}
\label{linear-regression}
\begin{array}{lcl}
&\min\limits_{\x\in\R^{n}}&\frac{1}{2}\|\x\|^2\\
&{\text{s.t.}}&\x\in\argmin\limits_{\z\in\R^{n}}\frac{1}{2}\left\|A\z-b\right\|^2.
\end{array}
\end{equation}

We compare the performance of our Bisec-BiO with several existing methods to solve this problem, namely averaging iteratively regularized gradient method (a-IRG) \citep{kaushik2021method}, bilevel gradient SAM method (BiG-SAM) \citep{sabach2017first}, minimal norm gradient method (MNG) \citep{beck2014first}, and dynamic barrier gradient descent method (DBGD) \citep{gong2021bi}. In this experiment, the feasible set of the lower-level problem is unbounded. Therefore, we cannot directly apply CG-BiO \citep{jiang2022conditional}. We use MATLAB function \texttt{lsqminnorm} to solve problem (\ref{linear-regression}) and obtain the optimal values $g^{*}$ and $p^{*}$ for benchmarking purposes.

\begin{figure}[!ht]
\begin{center}
\includegraphics[width=0.48\linewidth]{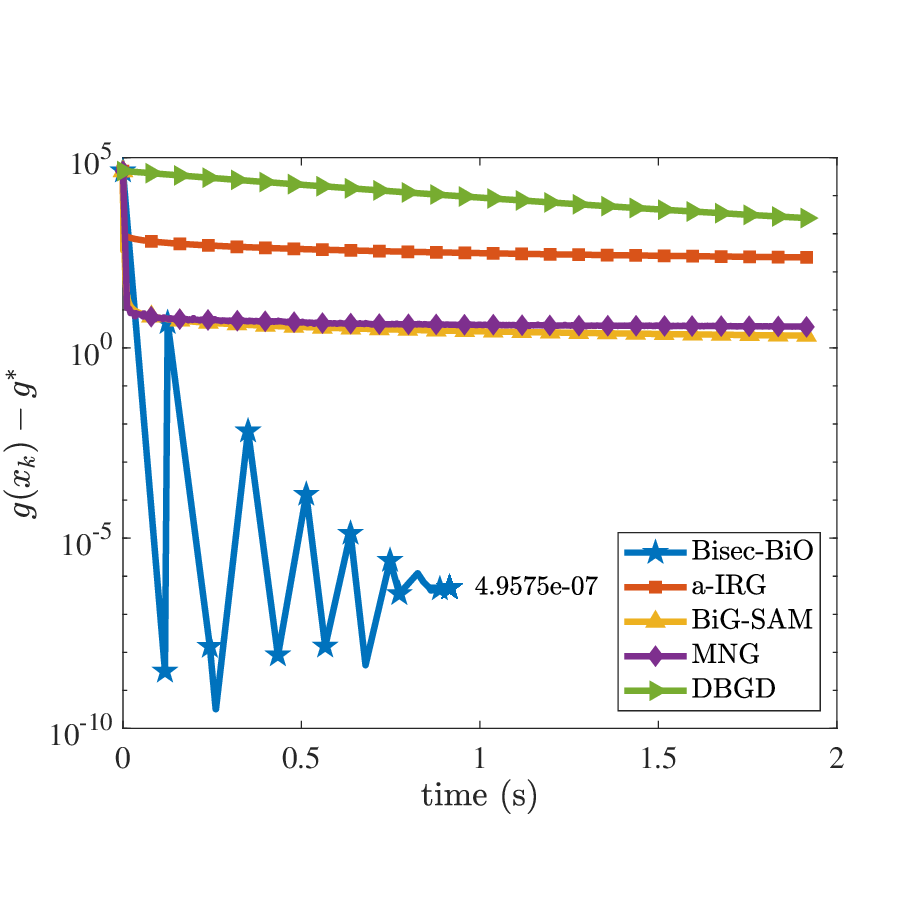}
\includegraphics[width=0.48\linewidth]{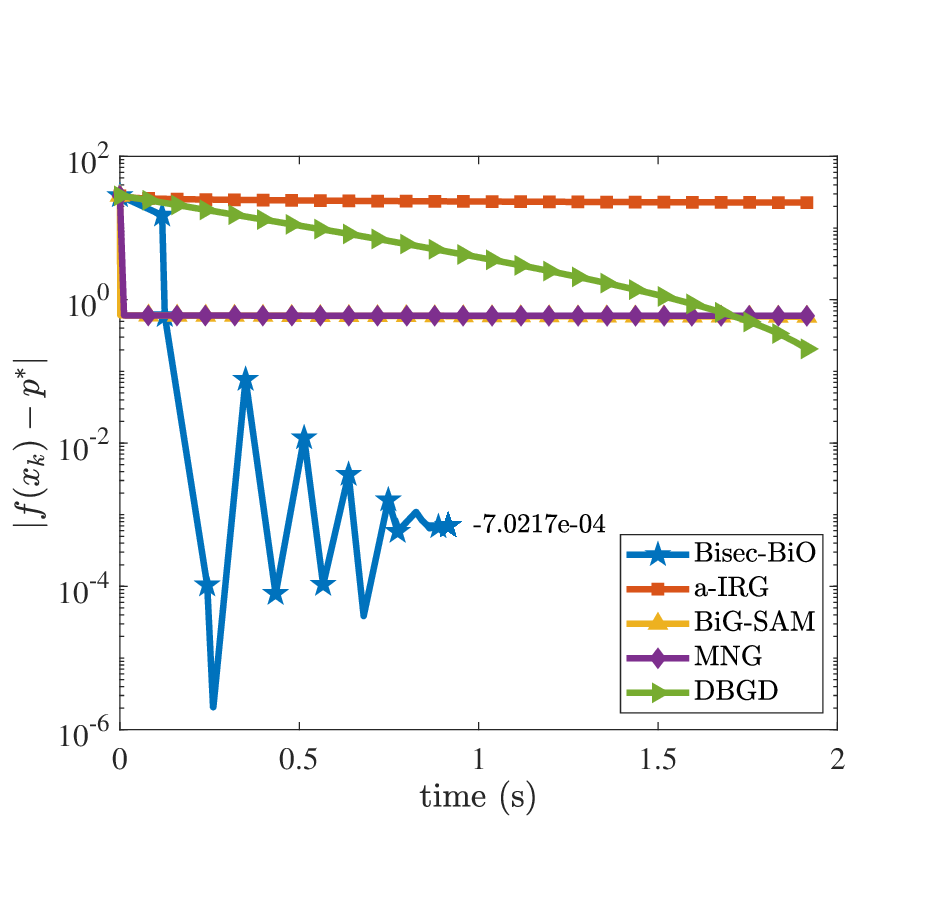}
\end{center}
\caption{The performance of Bisec-BiO compared with other methods in MNP.}
\label{linear-figure}
\end{figure}
As Figure \ref{linear-figure} demonstrates, our Bisec-BiO converges much faster than the other baseline methods for both the lower- and upper-level objectives, confirming our complexity results (see Table \ref{table1}). Our method is shown in the figure as an oscillating curve, which depends on the feasibility of System (\ref{system1}) at each iteration point and the update mode of our method. Moreover, the left and right subfigures in Figure \ref{linear-figure} show that the output of our algorithm meets the $(\epsilon_f, \epsilon_g)$-optimal solution criterion, as proven in Theorem \ref{thm1}.

\subsection{Logistic Regression Problem (LRP)}
\label{logistic regression problem}
We address the logistic regression binary classification problem using the 'a1a' dataset from LIBSVM\footnote{\url{https://www.csie.ntu.edu.tw/~cjlin/libsvmtools/datasets/binary/a1a.t}}, which contains $m=30,956$ instances, each with $n=123$ features. We randomly select a sample of $1,000$ instances and denote the feature matrix and labels as $A$ and $b$, respectively.

For the lower-level, we let $g_1(\x) = \frac{1}{m}\sum_{i=1}^{m}\log(1+\exp(-a_{i}^{\top}\x b_{i}))$, which has a $L_{g_1}$-Lipschitz continuous gradient with $L_{g_1} = \frac{1}{4m}\lambda_{\max}(A^TA)$. Here, $a_i$ represents an instance and $b_{i}\in\{-1,1\}$ is the corresponding label. Additionally, we set $g_2(\x)=I_C(\x)$, where $I_C(\x)$ is the indicator function of $C=\{\x\in \R^{n}: \|\x\|_1\le \lambda\}$ with $\lambda=10$. For the upper-level, we let $f_1(\x)=\frac{1}{2}\|\x\|^2$ and $f_2(\x)\equiv 0$, as well. We need to solve the following problem:
\begin{equation}
\label{logistic-regression}
\begin{array}{lcl}
&\min\limits_{\x\in\R^{n}}&\frac{1}{2}\|\x\|^2\\
&{\text{s.t.}}&\x\in\argmin\limits_{\z\in\R^{n}}\frac{1}{m}\sum_{i=1}^{m}\log(1+\exp(-a_{i}^{\top}\z b_{i}))\\
& &+I_C(\z).\\
\end{array}
\end{equation}

In this experiment, we compare the performance of our Bisec-BiO with the methods proposed in Section \ref{minimum norm solution problem} and the CG-based bilevel optimization method (CG-BiO) \citep{jiang2022conditional}. For benchmarking purposes, we utilize the Greedy FISTA algorithm \citep{liang2022improving} and MATLAB function \texttt{fmincon} to solve problem (\ref{linear-regression}) and obtain the optimal values $g^{*}$ and $p^{*}$, respectively. We employ the method proposed in \citet{liu2020projections} to compute the proximal operator of $h_c:=g_2+\delta_c$ in problem (\ref{pb:c}). Their method is demonstrated to have a worst-case complexity of $O(n^2)$ and an observed complexity of $O(n)$ in practice.

\begin{figure}[!ht]
\begin{center}
\includegraphics[width=0.48\linewidth]{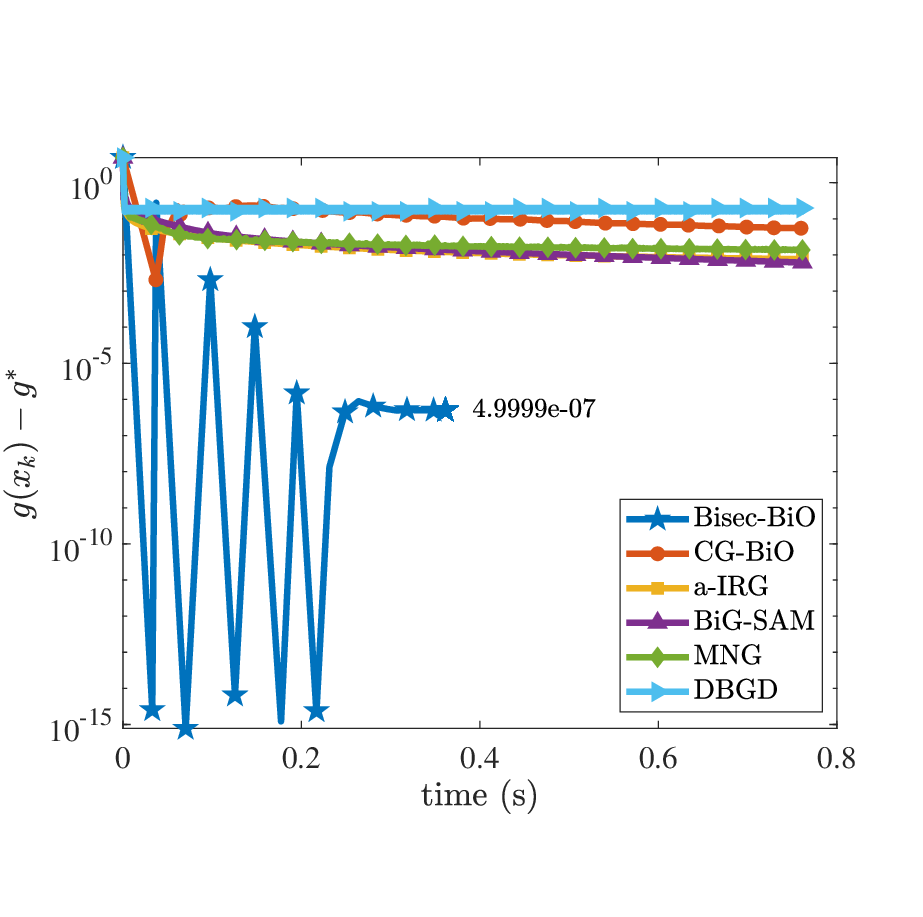}
\includegraphics[width=0.48\linewidth]{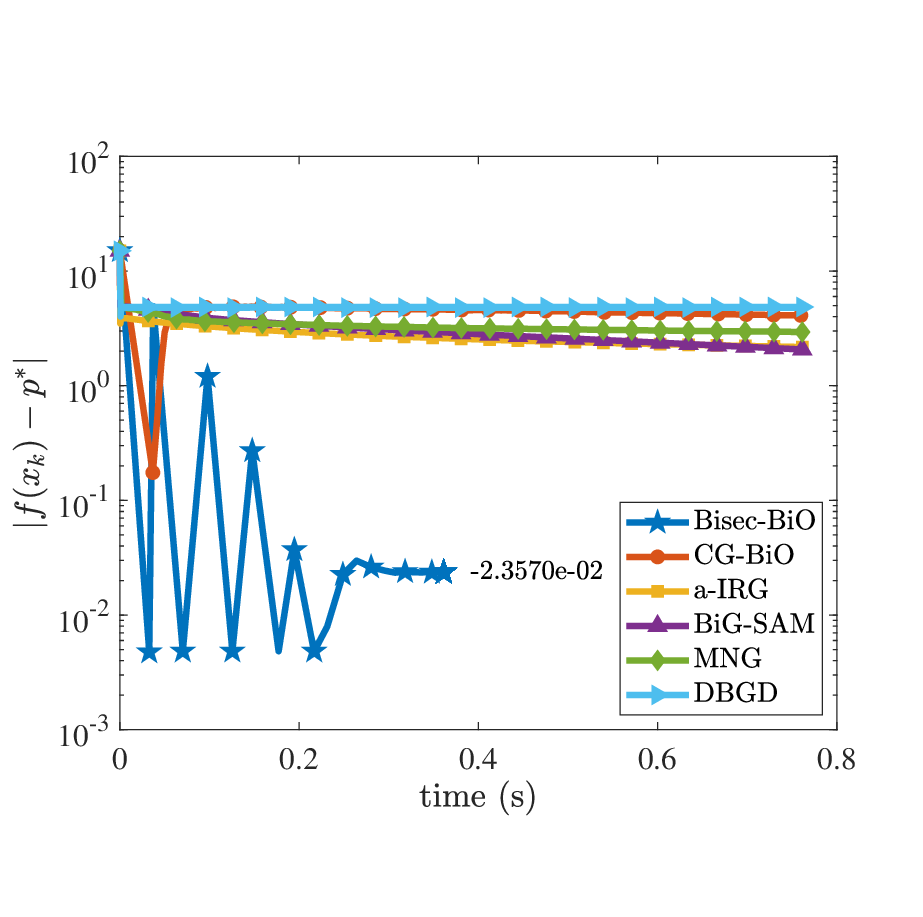}
\end{center}
\caption{The performance of Bisec-BiO compared with other methods in LRP.}
\label{logistic-figure}
\end{figure}
Figure \ref{logistic-figure} also shows that Bisec-BiO converges much faster than the other baseline methods for both the lower- and upper-level objectives. We have similar observations as in Figure \ref{linear-figure}.

\section{CONCLUSION}
\label{conclusion}
In this paper, we address the problem of minimizing a composite convex function at the upper-level within the optimal solution set of a composite convex lower-level problem. We introduce a bisection algorithm designed to discover an $\epsilon_f$-optimal solution for the upper-level objective and an $\epsilon_g$-optimal solution for the lower-level objective. Our method attains a near-optimal convergence rate of ${\tilde{\mathcal{O}}}\left(\max\{\sqrt{L_{f_1}/\epsilon_f},\sqrt{L_{g_1}/\epsilon_g}\}\right)$ for both upper- and lower-level objectives. Notably, this near-optimal rate aligns with the optimal rate observed in unconstrained smooth or composite optimization, neglecting the logarithmic term. We enhance convergence guarantees by imposing a H{\"{o}}lderian error bound assumption on the lower-level problem. Numerical experiments convincingly illustrate the substantial improvement our method offers over the state-of-the-art. In future research, we will explore the possibility of eliminating the logarithmic term from our complexity result.

\subsection*{Acknowledgements}
\label{Acknowledgements}
Rujun Jiang is partly supported by the National Key R\&D Program of China under grant 2023YFA1009300,
National Natural Science Foundation of China under grants 12171100 and 72394364, and Natural Science Foundation of Shanghai 22ZR1405100. Jiulin Wang is supported by China Postdoctoral Science Foundation under grants BX20220085 and 2022M710798.

\bibliographystyle{plainnat}
\bibliography{reference}

\vfill

\appendix

\onecolumn

\section{THE FAST ITERATIVE SHRINKAGE THRESHOLDING ALGORITHM (FISTA) AND CONVERGENCE RESULTS}
We adopt the Fast Iterative Shrinkage Thresholding Algorithm (FISTA) proposed in \citet{Beck2009} as the APG method (see Definition \ref{definition-APG}) to solve problem (\ref{pb:varphi}).
For any $L>0$, define $p_{L}(\y)$ as the unique minimizer of the following quadratic approximation of $\varphi(\x)$ at a given point $\y$:
\begin{equation*}
\label{Qxy}
p_{L}(\y):=\arg\min_{\x\in \R^n}  Q_L(\x,\y):=\left\{\varphi_1(\y)+\left\langle{\x-\y,\nabla {\varphi_1}(\y) }\right\rangle
+\frac{L}{2}\|\x-\y\|^2+\varphi_2(\x)\right\},
\end{equation*}
where $L$ actually plays the role of a step-size.

The classical FISTA scheme with a constant step-size for solving problem (\ref{pb:varphi}) is presented in Algorithm \ref{alg:fista}.
\begin{algorithm}
\caption{FISTA with constant step-size}\label{alg:fista}
\begin{algorithmic}[1]
\Require $L_{\varphi_1}$, $t_1 = 1$, $\y_{1}=\x_0 \in \R^n$
\For{\parbox[t]{0.12\linewidth}{$k=1,\cdots$}}
\State $\x_{k} = p_{L_{\varphi_1}}(\y_{k})$,
\State $t_{k+1} = \frac{1+\sqrt{1+4 t_{k}^{2}}}{2}$,
\State $\y_{k+1} = \x_{k} + \frac{t_{k}-1}{t_{k+1}}(\x_{k}-\x_{k-1})$,
\State $k=k+1$.
\EndFor
\end{algorithmic}
\end{algorithm}

As discussed in Section \ref{preliminaries}, a potential limitation of this classical scheme is its dependence on the knowledge or computation of the Lipschitz constant $L_{\varphi_1}$, which may not be practical. To overcome this issue, \cite{Beck2009} further proposed a variant of FISTA that incorporates a backtracking line search, we present it in Algorithm \ref{alg:fistalinesearch}.

\begin{algorithm}
\caption{FISTA with backtracking line search}\label{alg:fistalinesearch}
\begin{algorithmic}[1]
\Require $L_0>0$, $\eta>1$, $t_1 = 1$, $\y_{1}=\x_0 \in \R^n$
\For{\parbox[t]{0.12\linewidth}{$k=1,\cdots$}}
\State Find the smallest nonnegative integer value $i_k$ such that with ${\bar L} = \eta^{i_k}L_{k-1}$, 
\[
\varphi(p_{{\bar L}}(\y_{k}))\leq Q_{{\bar L}}(p_{{\bar L}}(\y_{k}),\y_{k}).
\]
\State $L_k = \eta^{i_k}L_{k-1}$,
\State $\x_{k} = p_{L_k}(\y_{k})$,
\State $t_{k+1} = \frac{1+\sqrt{1+4 t_{k}^{2}}}{2}$,
\State $\y_{k+1} = \x_{k} + \frac{t_{k}-1}{t_{k+1}}(\x_{k}-\x_{k-1})$,
\State $k=k+1$.

\EndFor
\end{algorithmic}
\end{algorithm}

The next lemma shows the convergence results of the objective function under the FISTA scheme (Algorithm \ref{alg:fista} or \ref{alg:fistalinesearch}).
\begin{lemma}[\citep{Beck2009}, Theorem 4.4.]\label{lemma}
Denote $X_{\varphi}^*$ as the optimal solution set of problem (\ref{pb:varphi}) and $\x_{\varphi}^*\in X_{\varphi}^*$ be any optimal solution. Let $\x_0^{\varphi}\in \R^n$ be an initial point. Let $\{\x_k\}$ be the sequence generated by the FISTA scheme (Algorithm \ref{alg:fista} or \ref{alg:fistalinesearch}). Then for any $k \ge 1$, we have
\begin{equation*}
\varphi(\x_k)- \varphi(\x_{\varphi}^*) \le \frac{2 \alpha L_{{\varphi}_1} }{(k + 1)^2} \|\x_0^{\varphi}-\x_{\varphi}^*\|^2,~\forall \x_{\varphi}^*\in X_{\varphi}^*.
\end{equation*}
\end{lemma}
Here, $\alpha=1$ for the classical FISTA scheme (Algorithm \ref{alg:fista}), and $\alpha=\eta$ for the FISTA scheme with a backtracking line search (Algorithm \ref{alg:fistalinesearch}), where $\eta$ is the backtracking parameter. This lemma demonstrates that an $\epsilon$-optimal solution of problem (\ref{pb:varphi}) can be obtained by Algorithm \ref{alg:fista} or \ref{alg:fistalinesearch} within at most $\mathcal{O}(\sqrt{L_{\varphi_1}/{\epsilon}})$ iterations.

\section{MOTIVATING EXAMPLES}
Many applications in machine learning and signal processing involve regularized problems, where the upper-level objectives represent the regularization terms, and the lower-level objectives consist of the loss functions and the additional constraint terms. We present some motivating examples below.
\begin{example}[Minimum Norm Solution of Least Squares Regression Problem (MNP)]
\label{ex1}
Linear inverse problems aim to reconstruct a vector $\x\in\mathbb{R}^n$ from a set of measurements $b\in\mathbb{R}^m$ that satisfy the following relation: $b = A\x + \rho \varepsilon$, where $A:\mathbb{R}^n\rightarrow\mathbb{R}^m$ is a given linear mapping, $\varepsilon\in\mathbb{R}^m$ denotes an unknown noise vector, and $\rho>0$ denotes its magnitude. Linear inverse problems can be solved using various optimization techniques, and we focus on the bilevel formulation \citep{beck2014first,sabach2017first,dempe2021simple,latafat2023adabim,merchav2023convex}. 

The lower-level objective function in this formulation is given by 
\begin{equation*}
\begin{split}
g(\x) = \frac{1}{2}\left\|A\x-b\right\|^2 + I_C(\x),
\end{split}
\end{equation*}
where the set $C$ is a closed convex set that can be chosen as $C=\mathbb{R}^{n}$, $C=\{\x\in \R^n: \x\ge 0\}$ or $C=\{\x\in \R^n: \|\x\|_1\le \lambda\}$ for some $\lambda>0$,
and $I_C(\x)$ is the indicator function of the set $C$, defined as $I_C(\x)=0$ if $\x \in C$ and $I_C(\x)=+\infty$ if $\x \notin C$.

This problem may have multiple optimal solutions. Therefore, a natural choice is to consider the minimal norm solution problem, which seeks to find the optimal solution with the smallest Euclidean norm \citep{beck2014first,sabach2017first,latafat2023adabim}:
\begin{equation*}
\begin{split}
f(\x) = \frac{1}{2}\left\|\x\right\|^2.
\end{split}
\end{equation*}
We then solve the bilevel optimization problem:
\begin{equation*}
\begin{array}{lcl}
&\min\limits_{\x\in\mathbb{R}^{n}}&\frac{1}{2}\|\x\|^2\\
&{\rm{s.t.}}&\x\in\argmin\limits_{\z\in\R^{n}}\frac{1}{2}\left\|A\z-b\right\|^2 + I_C(\z).
\end{array}
\end{equation*}
\end{example}
For this example, the proximal mapping of $h_c$ reduces to an orthogonal projection onto the $\ell_2$-norm ball when $C = \mathbb{R}^n$. This scenario corresponds to the experiment in Section \ref{minimum norm solution problem}.

When $C=\{\x\in \R^n: \x\ge 0\}$, we have
\[
{\rm prox}_{h_c}(\y)=\frac{\sqrt{2c}}{\max\{\|P_C(\y)\|,\sqrt{2c}\}}\cdot P_C(\y),
\]
where $P_C(\y)=\max(\y,0)$. This result
is an implementation of Theorem 7.1 in \cite{bauschke2018projecting}.

When $C=\{\x\in \R^n: \|\x\|_1\le \lambda\}$, the proximal mapping of $h_c$ simplifies to an orthogonal projection onto the intersection of a $\ell_2$-norm ball and a $\ell_1$-norm ball. This projection has a worst-case complexity of $O(n^2)$ and an observed complexity of $O(n)$ in practice \citep{liu2020projections}. In this case, the lower-level objectives satisfy the H{\"{o}}lderian error bound assumption (Assumption \ref{as3}) with $r=2$.

\begin{example}[Sparse Solution of Least Squares Regression Problem (SSP)]
\label{ex1-1}
Considering the same settings as in Example \ref{ex1}. To simplify the model and save computational resources and efficiency, we seek to reduce the number of features in the vector $\x\in\mathbb{R}^n$ that minimizes the linear inverse regression function $g(\cdot)$. This means that our goal is to find a sparse solution among all the minimizers of $g(\cdot)$. Therefore, any function that promotes sparsity can be used for this purpose. For example, the well-known elastic net regularization is a good choice \citep{zou2005regularization,friedlander2008exact,de2009elastic,rodola2013elastic,amini2019iterative,merchav2023convex}. The elastic net regularization is defined as
\begin{equation*}
\begin{split}
f(\x) = \left\|\x\right\|_1 + \frac{\alpha}{2}\left\|\x\right\|^2,
\end{split}
\end{equation*}
where $\alpha>0$ regulates the trade-off between $\ell_1$ and $\ell_2$ norms.
\end{example}
For this example, the proximal mapping of $h_c$ reduces to an orthogonal projection onto the elastic-net constraints when $C = \mathbb{R}^n$ \citep{duchi2009elastic,mairal2010online,gong2011efficient,rodola2013elastic}.

\begin{example}[Logistic Regression Classification Problem (LRP)]
\label{ex2}
The goal of a binary classification problem is to establish a mapping from the feature vectors $a_i$ to the target labels $b_i$. A common machine learning approach for this task is to minimize the logistic regression function of the given dataset \citep{amini2019iterative,gong2021bi,jiang2022conditional,latafat2023adabim,merchav2023convex}. More precisely, we have a feature matrix $A \in \mathbb{R}^{m \times n}$ and a corresponding label vector $b \in \mathbb{R}^{m}$, where each $b_i \in {-1, 1}$. The logistic loss function is then given by
\begin{equation*}
\begin{split}
g_1(\x) = \frac{1}{m}\sum_{i=1}^{m}\log(1+\exp(-a_{i}^{\top}\x b_{i})).
\end{split}
\end{equation*}
Over-fitting may occur when the number of features is not negligible relative to the number of instances $m$. A common way to address this problem is to regularize the logistic objective function with a specific function or add a constraint \citep{jiang2022conditional,merchav2023convex}. For example, we can choose $g_2(\x) = I_C(\x)$, where $I_C(\x)$ is the indicator of the set $C=\{\x\in \R^n: \|\x\|_1\le \lambda\}$ as shown in Example \ref{ex1}.

Similarly, this problem may have multiple optimal solutions. Therefore, it is natural to consider the minimal norm solution problem with the smallest Euclidean norm as described in Example \ref{ex1}. Now, we need to solve the following bilevel optimization problem \citep{gong2021bi,jiang2022conditional,latafat2023adabim}:
\begin{equation*}
\begin{array}{lcl}
&\min\limits_{\x\in\mathbb{R}^{n}}&\frac{1}{2}\|\x\|^2\\
&{\rm{s.t.}}&\x\in\argmin\limits_{\z\in\R^{n}}\frac{1}{m}\sum\limits_{i=1}^{m}\log(1+\exp(-a_{i}^{\top}\z b_{i})) + I_C(\z).
\end{array}
\end{equation*}
\end{example}
The proximal mappings of $h_c$ are the same for the choices of the set $C$ in Example \ref{ex1}. In particular, when $C=\{\x\in \R^n: \|\x\|_1\le \lambda\}$, this scenario corresponds to the experiment in Section \ref{logistic regression problem}.

\begin{example}[Sparse Solution of Logistic Regression Classification Problem (SSLRP)]
\label{ex2-1}
Considering the same settings as in
Example \ref{ex2}, and we seek to reduce the features in the vector $\x\in\mathbb{R}^n$ that minimizes the logistic regression function with a regularization term. For this purpose, we can choose the elastic-net regularization, which is proposed in Example \ref{ex1-1}.
\end{example}
Similarly, for this example, the proximal mapping of $h_c$ is an orthogonal projection onto the elastic-net constraints when $C = \mathbb{R}^n$, as shown in Example \ref{ex1-1}.

\section{PROOF OF THE MAIN THEOREMS}
\subsection{Proof of Theorem \ref{thm1}}
\begin{proof} 
We first show that ${\hat \x}$ is an $(\epsilon_f,\epsilon_g)$-optimal solution of problem (\ref{pb:BP}).
In Step 10, we let $c=u$. Then Condition (\ref{approxpsi}) is not satisfied. According to Lemma \ref{lem:system}, the inequality $g({\hat \x})\le g^*+\epsilon_g$ holds. Next, we prove $f({\hat \x})\le p^*+\epsilon_f$ also holds. We break the proof into two cases. 

Case (1), if $u\le p^*$, then we have $f({\hat \x}) \le u\le p^*$.

Case (2), if $u> p^*$, then $p^*$ lies on $[l,u]$ since $l\le p^*$ always holds. Therefore, we have
\[
f({\hat \x}) \le u\le u+p^*-l\le p^*+\epsilon_f,
\]
where the last inequality is from the stop criterion that $u-l\le \epsilon_f$. To sum up, the point ${\hat \x}$ is an $(\epsilon_f,\epsilon_g)$-optimal solution of problem (\ref{pb:BP}). In the following, we present the total operation complexity of our method.

In Step 1, we invoke APG oracles to obtain initial bounds $l$ and $u$. To obtain the initial lower bound of $p^*$, we invoke the APG oracle ${\tilde{\x}_f}$=APG$(f_1, f_2, L_{f_1}, \x_0^f, \epsilon_f/2)$ to solve problem (\ref{pb:fvalue}). By Lemma \ref{lemma}, this can be done within ${\mathcal{O}}\left(\sqrt{{ L_{f_1}}/{ \epsilon_f}}\right)$ iterations. The corresponding initial lower bound is $l=f(\tilde{\x}_{f})-\epsilon_f/2$. To obtain the initial upper bound of $p_{\epsilon_g}^*$, we invoke the APG oracle ${\tilde{\x}_g}$=APG$(g_1, g_2, L_{g_1}, \x_0^g, \epsilon_g/2)$ to solve problem (\ref{pb:g}). Similarly, we can approximately solve problem (\ref{pb:g}) within
${\mathcal{O}}\left(\sqrt{ {L_{g_1}}/{ \epsilon_g}}\right)$ iterations. The corresponding initial upper bound is given by $u=f(\tilde{\x}_{g})$. 

In Step 3, we invoke APG oracle ${\tilde{\x}_c}$=APG$(g_1, h_c, L_{g_1}, \x_0^c, \epsilon_g/2)$ to
solve problem (\ref{pb:system}). According to Lemma \ref{lemma}, this can be done within ${\mathcal{O}}\left(\sqrt{ {L_{g_1}}/{\epsilon_g}}\right)$ iterations. The number of invoking APG oracle does not exceed
\[
\left\lceil\log_2\frac{u-l}{\epsilon_f} \right\rceil_{+}=\left\lceil\log_2\frac{f(\tilde{\x}_{g})-f(\tilde{\x}_{f})+\epsilon_f/2}{\epsilon_f} \right\rceil_{+}={\mathcal{O}}\left(\log \frac{1}{\epsilon_f}\right)
\]
where $l$ and $u$ are initial lower and upper bounds, and $\lceil a \rceil_{+}$ represents the smallest non-negative integer that is no less than $a$.

Thus, the total number of evaluations of the function values $f_1$, $f_2$, $g_1$, and $g_2$, the gradients $\nabla f_1$ and $\nabla g_1$, and the calls of the proximal mapping concerning function $h_c$ does not exceed $T$, where
\[
T={\mathcal{O}}\left(\sqrt{\frac {L_{f_1}}{ \epsilon_f}}\right)+{\mathcal{O}}\left(\sqrt{\frac {L_{g_1}}{ \epsilon_g}}\right)+{\mathcal{O}}\left(\sqrt{\frac {L_{g_1}}{ \epsilon_g}}\right)\cdot {\mathcal{O}}\left(\log \frac{1}{\epsilon_f}\right)=\tilde {\mathcal{O}}\left(\max\left\{\sqrt{\frac {L_{f_1}}{ \epsilon_f}},\sqrt{\frac {L_{g_1}}{ \epsilon_g}} \right\}\right),
\]
where $\tilde {\mathcal{O}}$ suppresses a logarithmic term.
\end{proof}

\subsection{Proof of Proposition \ref{prop1}}
\begin{proof}
By Assumption \ref{as2}(i), the set $X_g^*$ is closed and compact. Then we can let $\hat{\x}^* = \argmin\limits_{\x \in X_g^*} \|\x - \hat{\x}\|$ such that $\|\hat{\x}- \hat{\x}^* \| = {\rm dist}(\hat{\x}, X_g^*)$. It can be easily demonstrated that $X_g^*$ is a convex set, ensuring the well-definedness of $\hat{\x}^*$. By Assumption \ref{as3}, we have
\begin{equation}\label{boundx}
\frac{\alpha}{r} \|\hat{\x}- \hat{\x}^* \|^r \leq g(\hat{\x}) - g^* \leq \epsilon_g \implies \|\hat{\x}- \hat{\x}^* \| \leq \left(\frac{r\epsilon_g}{\alpha}\right)^{\frac{1}{r}}.    
\end{equation}

By Assumption \ref{as2}, it follows that $f$ is $B_f$-Lipschitz continuous on ${\rm dom}(g_2)$ (see (\ref{Lip})). Combining this result with (\ref{boundx}), we have
\[
f(\hat {\x})-p^*\ge f(\hat {\x})-f(\hat{\x}^*) \geq -B_f\|\hat{\x}- \hat{\x}^* \|
\ge -B_f \left(\frac{r\epsilon_g}{\alpha}\right)^{\frac{1}{r}},
\]
where the first inequality is from that $p^*$ is the optimal value of problem \eqref{pb:BP} and $\hat{\x}^*\in X_g^*$.
\end{proof}

\subsection{Proof of Corollary \ref{cor1}}
\begin{proof}
Let $\hat {\x}$ be the output of Algorithm \ref{main algorithm}. Then $\hat {\x}$ is an $(\epsilon_f,\epsilon_g)$-optimal solution of (\ref{pb:BP}) satisfying
\[
f(\hat{\x}) - p^*\le \epsilon_f, ~ g(\hat{\x}) - g^*\le \epsilon_g.
\]
By Proposition \ref{prop1} and the setting $\epsilon_g=\frac{\alpha}{\gamma} \left(\frac{\epsilon_f}{B_f}\right)^r$, we have
\[
f(\hat {\x})-p^*
\ge -B_f \left(\frac{r\epsilon_g}{\alpha}\right)^{\frac{1}{r}}=-\epsilon_f.
\] 
Then according to Theorem \ref{thm1}, with an operation complexity $\tilde {\mathcal{O}}\left({1/{\sqrt{\epsilon_f^r} }}\right)$, we can find an $\hat {\x}$ such that
\[
|f(\hat {\x})-p^*|\le \epsilon_f,~ g(\hat {\x})-g^*\le \epsilon_g.
\]
\end{proof}

\subsection{Proof of Proposition \ref{prop2}}
\begin{proof}
Since $\hat{\x}$ satisfies $g(\hat{\x})-g^*\le \epsilon_g$, then $\hat{\x}$ is a feasible solution of (\ref{pb:BP2}) with $\epsilon=\epsilon_g$, and $f(\hat{\x})$ is an upper bound of $p_{\epsilon_g}^*$ (see Lemma \ref{lem:system}). Applying Assumption \ref{as4} with $\epsilon=\epsilon_g$, we have
\[
f(\hat {\x})-p^*\ge p_{\epsilon_g}^*-p^*=v(\epsilon_g)-v(0)\ge -L_{\epsilon_g}\epsilon_g.
\]
\end{proof}

\subsection{Proof of Corollary \ref{cor2}}
\begin{proof}
This proof closely resembles the one presented in Corollary \ref{cor1}. Nevertheless, for the sake of thoroughness, we provide a complete exposition. Let $\hat {\x}$ be the output of Algorithm \ref{main algorithm}. Then $\hat {\x}$ is an $(\epsilon_f,\epsilon_g)$-optimal solution of (\ref{pb:BP}) satisfying
\[
f(\hat{\x}) - p^*\le \epsilon_f, ~ g(\hat{\x}) - g^*\le \epsilon_g.
\]
By Proposition \ref{prop2} and the setting $\epsilon_f=L_{\epsilon_g}\epsilon_g$, we have
\[
f(\hat {\x})-p^*
\ge -L_{\epsilon_g}\epsilon_g=-\epsilon_f.
\] 
Then according to Theorem \ref{thm1}, with an operation complexity $\tilde {\mathcal{O}}\left( \max \{1/\sqrt{L_{\epsilon_g}\epsilon_g}, 1/\sqrt{\epsilon_g}\}\right)$, we can find an $\hat {\x}$ such that
\[
|f(\hat {\x})-p^*|\le \epsilon_f,~ g(\hat {\x})-g^*\le \epsilon_g.
\]
\end{proof}

\section{ADDITIONAL NUMERICAL EXPERIMENTS}
\label{linear-elastic}
We consider the SSP problem from Example \ref{ex1-1} on the YearPredictionMSD dataset, setting $C = \mathbb{R}^n$ and $\alpha = 0.02$. Our objective is to solve the following bilevel problem:
\begin{equation}
\label{linear-regression-elastic}
\begin{array}{lcl}
&\min\limits_{\x\in\R^{n}}& \frac{\alpha}{2}\|\x\|^2 + \|\x\|_1\\
&{\text{s.t.}}&\x\in\argmin\limits_{\z\in\R^{n}}\frac{1}{2}\left\|A\z-b\right\|^2.
\end{array}
\end{equation}

We compare the performance of our Bisec-BiO method with the averaging iteratively regularized gradient (a-IRG) method \citep{kaushik2021method}. It's worth noting that a-IRG can handle non-smooth upper-level objectives, a capability that other methods lack. For benchmarking purposes, we use the MATLAB functions \texttt{lsqminnorm} and \texttt{fmincon} to obtain the optimal values $g^{\ast}$ and $p^{\ast}$, respectively. Moreover, we adopt the method proposed by \citet{gong2011efficient} to compute the proximal operator of $h_c:=g_2+\delta_c$ in problem (\ref{pb:c}). The other settings are consistent with Section \ref{minimum norm solution problem}.
\begin{figure}[!ht]
\begin{center}
\includegraphics[width=0.25\linewidth]{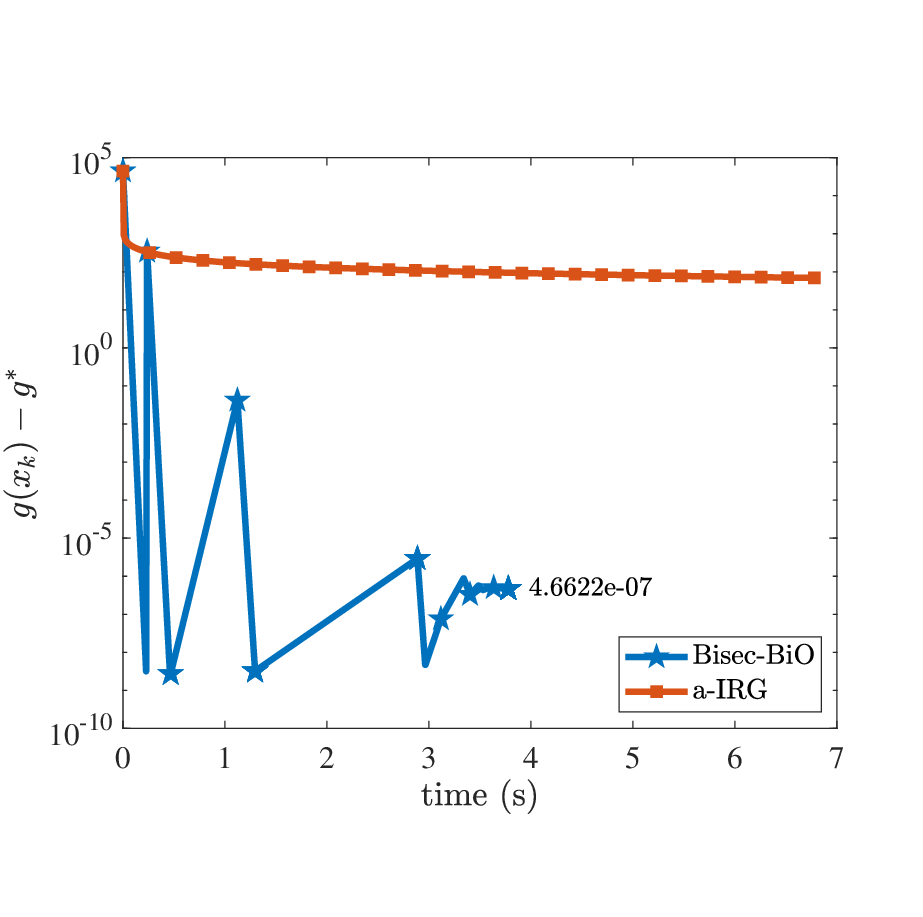}
\includegraphics[width=0.25\linewidth]{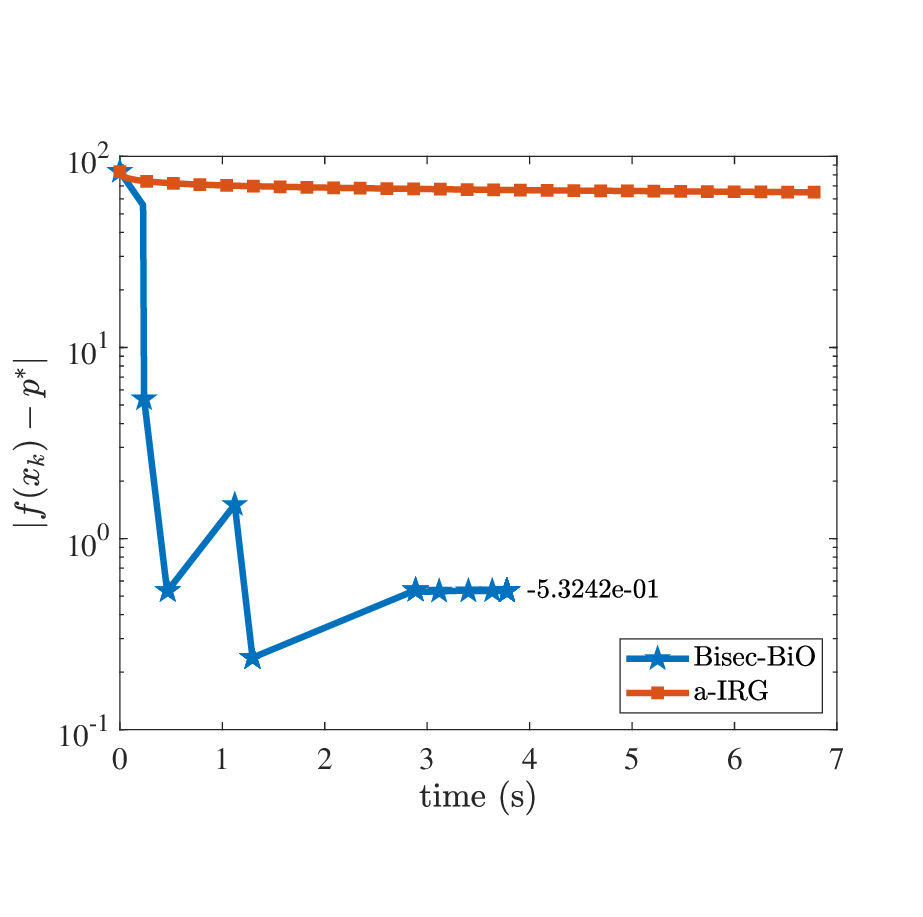}
\end{center}
\caption{The performance of Bisec-BiO compared with other methods in SSP.}
\label{linear-elastic-figure}
\end{figure}

As illustrated in Figure \ref{linear-elastic-figure}, Bisec-BiO converges significantly faster than a-IRG for both lower- and upper-level objectives. The numbers on the right-hand side of the last iterates of our method (denoted as $\hat{\x}$) represent the differences between $g(\hat{\x}) - g^*$ and $f(\hat{\x})- p^*$, respectively. This confirms that $\hat{\mathbf{x}}$ is an $(\epsilon_f, \epsilon_g)$-optimal solution.

\end{document}